\documentclass[12pt]{article}

\usepackage[a4paper, left=3.2cm,top=3.2cm]{geometry}
\usepackage{authblk}
\usepackage[onehalfspacing]{setspace}

\usepackage{amsmath,amssymb,amsthm}
\usepackage{booktabs}
\usepackage{mathtools}
\usepackage{enumitem}

\usepackage[colorlinks=true,bookmarks=false,citecolor=blue,urlcolor=blue]{hyperref}

\usepackage{siunitx}

\DeclareMathOperator{\Wo}{\mathbf{W}}
\DeclareMathOperator{\Uo}{\mathbf{U}}
\DeclareMathOperator{\Fo}{\mathbf{F}}

\DeclareMathOperator{\So}{\mathbf{S}}
\DeclareMathOperator{\Ao}{\mathbf{A}}

\newcommand{\R}{\mathbb{R}}
\newcommand{\Z}{\mathbb{Z}}
\newcommand{\edot}{\,\cdot \,}
\newcommand{\BB}{\mathcal{B}}

\newcommand{\rpos}{\mathbf{r}}
\newcommand{\spos}{\mathbf{s}}

\newcommand{\x}{\mathbf{x}}

\DeclarePairedDelimiter{\abs}{\lvert}{\rvert}
\DeclarePairedDelimiter{\norm}{\lVert}{\rVert}

\newtheorem{theorem}{Theorem}

\newtheorem{definition}[theorem]{Definition}

\allowdisplaybreaks

\title{Sampling and resolution in sparse view photoacoustic tomography}

\date{}

\author{Markus Haltmeier}
\author{Daniel Obmann}

\affil{Department of Mathematics, University of Innsbruck\authorcr
Technikerstrasse 13, 6020 Innsbruck, Austria
 \authorcr E-mail:  \texttt{markus.haltmeier@uibk.ac.at}
 }

\author{Karoline Felbermayer}
\author{Florian Hinterleitner}
\author{Peter Burgholzer}

\affil{Research Center for Non-Destructive Testing (RECENDT)\authorcr Altenberger Stra{\ss}e 69, 4040 Linz, Austria
 \authorcr E-mail:  \texttt{peter.burgholzer@recendt.ac.at}
 }

\begin{document}

\maketitle

\begin{abstract}
We investigate resolution in photoacoustic tomography (PAT). Using Shannon theory, we investigate the theoretical resolution limit of sparse view PAT theoretically, and empirically demonstrate that all reconstruction methods used exceed this limit.
\end{abstract}

\section{Introduction}

The resolution and accuracy of photoacoustic tomography (PAT) depends on various factors including acoustic attenuation, limited bandwidth of the detection system and the number of available data samples.  In this work we investigate sampling and resolution of PAT from angularly undersampled data.  We analyze the theoretically achievable resolution given a maximal bandwidth $\Omega$ of the data. We derive conditions how to sample in the temporal and angular direction for 2D PAT using a circular arrangement of sensors. In sparse view PAT, the temporal sampling condition is met, while in the angular direction data are undersampled. As a consequence, not all objects in the class of  functions with  bandwidth $\Omega$  can be recovered and undersampling artefacts are introduced.

\begin{figure}[htb!]
\centering
\includegraphics[width=0.8\textwidth]{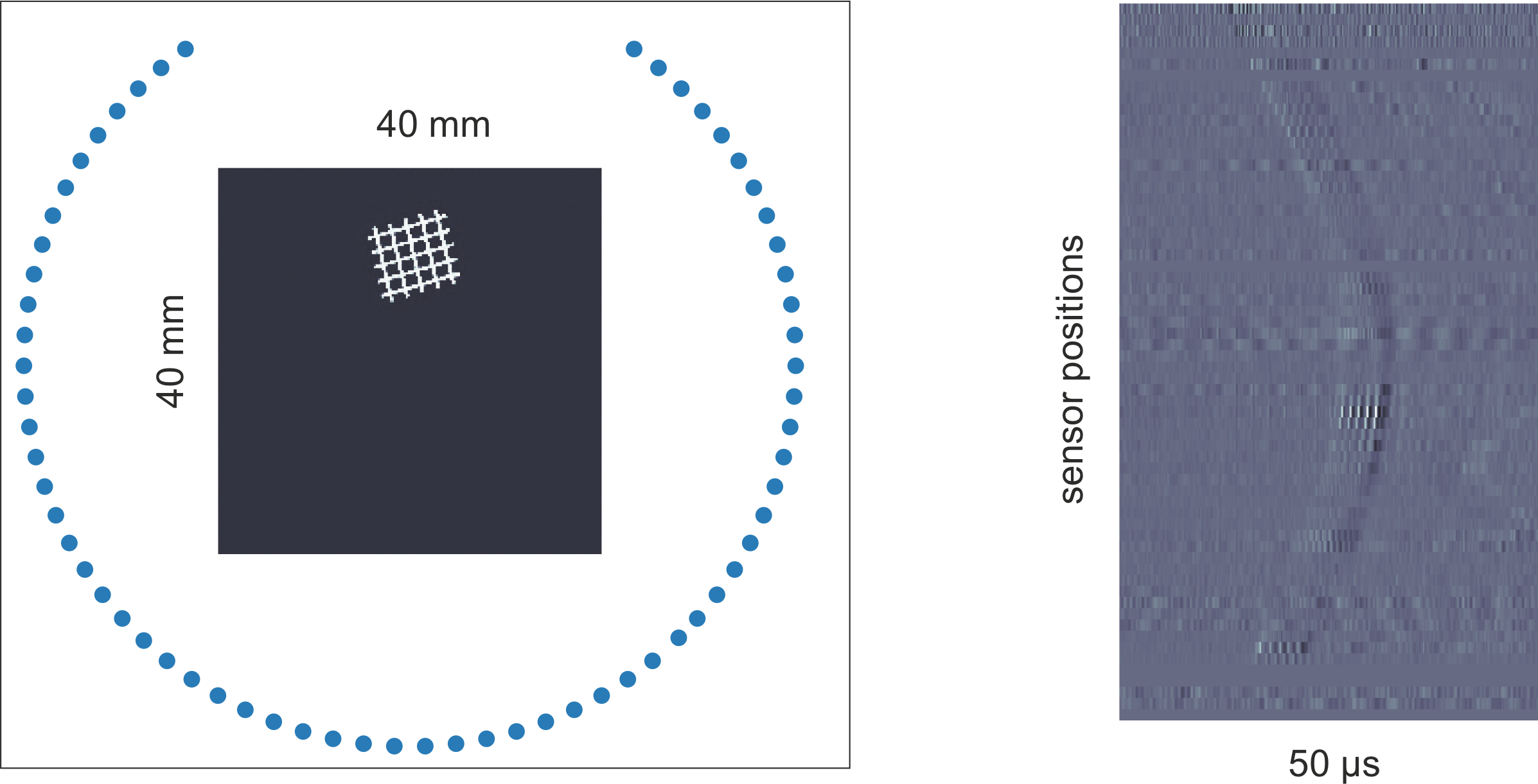}
\caption{Left: PAT data are collected with 64 sensors located on a circle of radius \SI{40}{\mm}. The initial pressure of the grid like phantom is contained in a $40 \times  40$ \si{\mm\squared} square centered at the midpoint of the circle. Right: The detected signal amplitude of the 64 detector points over the time $50 \si{\micro\second}$.} \label{fig:setup}
\end{figure}

In any application of PAT, the class objects to be reconstructed is not an arbitrary object of bandwidth $\Omega$. Instead it obeys additional structure and regularity that may or may not be available explicitly. In such a situation, the resolution can be significantly higher than indicated by the angular sampling condition. This reflects common practice in PAT that angular undersampling is used. In particular, using nonlinear reconstruction methods  the reconstruction quality significantly depends on the class of objects to be reconstructed.  In this paper  we consider a particular class of objects with a grid-like structure having well-definable resolution, see Figure~\ref{fig:setup}. Via numerical simulations we investigate if sparse view PAT data is capable to resolve the signal class. We compare standard quadratic Tikhonov regularization without specific prior, joint $\ell^1$ regularization using  sparsity and  positivity prior, and deep learning based reconstruction methods using training data as prior. We demonstrate that despite the angular undersampling, all methods are capable of well resolving the grid-like structure.

\section{Theory}

We consider  2D PAT imaging model in circular geometry. Let $f \colon \R^2\to\R$ denote  the PA source (initial pressure distribution). The induced  pressure wave  satisfies the wave equation $ \partial^2p(\rpos, t) -  c^2 \Delta_\rpos p(\rpos, t) = \delta^\prime(t) f(\rpos) \quad \text{for } (\rpos, t) \in \R^2\times \R$, where $\rpos\in\R^2$ is the spatial location, $t \in \R$ the  time, $\Delta_\rpos$ the spatial Laplacian, and $c$ is the constant speed of sound. After rescaling time $ t \gets ct$ we assume  $c=1$ in the following. We assume $p(\rpos,t)=0$ for $t<0$ such  that the solution $p(\rpos,t)$ is uniquely defined and denoted by $\Wo f$. We  assume  that the  acoustic pressure is  measured with point like  sensors on  the circle  $S_R = \{ \x \in \R^2 \mid \norm{\x} = R\}$ with radius $R>0$, each having spatial impulse response function (IRF)  $\varphi_\Omega \colon \R \to \R$, where $\Omega$ is the essential bandwidth determining the resolution.  For the experimental study we limit the bandwidth by convolving  data  with a Gaussian filter. The aim is to recover $f$ from samples of $\phi_\Omega  \ast_t \Wo f $ made on $S_R$.

\subsection{Resolution}

The IRF fully guides the achievable spatial resolution in PAT.
For theoretical analysis we assume that $\phi_\Omega \in L^2(\R) \cap L^1 (\R)$ is an even function with $\norm{\Fo_t  \phi_\Omega}_\infty = 1$.  Moreover we define $\Phi_\Omega$ by $\Fo \Phi_\Omega = \Fo_t  \phi_\Omega$  and refer to it as the points spread function (PSF). Here $\Fo$ and  $\Fo_t$ denote the Fourier transform in the spatial  and temporal variable, respectively.    In \cite{zangerl2020multi} the convolution identity $\phi_\Omega \ast \Wo f =  \Wo [\Phi_\Omega \ast f ]$ has been derived  relating the IRF and the PSF.

\begin{definition}\label{thm:res}
Let  $a >0$. A  subspace  $U \subseteq L^2(\R^2)$  is  called $a$-resolved by   $\Phi_\Omega$, if  $ \norm{  \Phi_\Omega \ast f }^2 \geq a \norm{  f }^2 $ for all $f \in U$. Likewise a subspaces  $V \subseteq L^2(S_R \times \R)$ is called $a$-resolved by  $\phi_\Omega$, if $ \norm{  \phi_\Omega \ast_t  g }^2 \geq   a \norm{  g }^2 $ for all $g \in V$. 
\end{definition} 

\begin{theorem}[Resolution] \label{thm:res}
Subspace $U$ is $a$-resolved by $\Phi_\Omega$ if and only if  $V = \Wo(U)$ is $a$-resolved by $\phi_\Omega$.
\end{theorem}

The relevance of  Theorem \ref{thm:res} for the resolution is most easily illustrated for the ideal low pass filter where $ \Fo_t \phi_\Omega (\omega) = 1$ for $\abs{\omega} \leq \Omega$ and zero otherwise.  Then $ \Fo \Phi_\Omega (\xi) = 1$ for $\norm{\xi} \leq \Omega$ and zero otherwise.
The largest space that is resolved by $\Phi_\Omega$ is the space $\BB_\Omega(\R^2)$ of $\Omega$ band-limited functions. Theorem \ref{thm:res}  states that  $\Wo(\BB_\Omega(\R^2))$ is resolved by $\phi_\Omega$. Likewise $\Wo( \BB_{b}(\R^2) )$ is not resolved by $\phi_\Omega$ if $b > \Omega$, hence $\Omega$ exactly characterizes the spatial resolution induced by the  ideal low pass.

\subsection{Sampling}

The derivation of sampling conditions requires fixing a space of functions where sampling is applied. Here we work with the space of band-limited functions and equidistant sampling $h_t>0$ in time.  Suppose that $\phi_\Omega$ resolves the space of $\BB_\Omega(\R)$ of band-limited functions, which implies that the same holds for $\Phi_\Omega$. We will sloppily say that $\Omega$ is the bandwidth of $\phi_\Omega$. Let $(u_k)_{k \in \Z^2}$ be a frame of $\BB_\Omega (\R^2)$, $\Uo^\ast$ be the synthesis operator and $h_t > 0$ the temporal step size. We define the spatially and temporally discretized operator 
$    \Wo_{\Omega, h_t} \colon \ell^2( \Z^2 )  \to L^2( S_R \times  \Z )$ by $   \Wo_{\Omega, h_t} (  x )  =  ( (\Wo \circ \Uo^\ast)(x)( \edot , m h_t) )_{m \in \Z} $. The basic question of sampling theory is finding conditions on the step size $h_t$ such that $\Wo_{\Omega, h_t} x $ uniquely determines $ \Uo^\ast(x) \in  \BB_\Omega (\R^2)$.

\begin{theorem}[Temporal sampling] \label{thm:temporal}
$\Wo_{\Omega, h_t} (x)$ uniquely determines $\Uo^\ast(x)$ if and only  if $h_t \leq \pi / \Omega$.
\end{theorem}

For the  spatial sampling conditions we take the frame $(u_k)_{k \in \Z^2}$ formed by  the translates $u_k(x) = \Phi_\Omega( x - k h_x )$ of the low pass filter $\Phi_\Omega$. The multi-dimensional Shannon sampling theorem gives the following. 

\begin{theorem}[Spatial sampling] \label{thm:spatial}
Any $f \in \BB_\Omega (\R^2)$ is  the form $f = \Uo^\ast(x)$ if and only  if $h_x \leq \pi / \Omega$.
\end{theorem}

Angular sampling crucially depends on the  location of the function to be recovered. 
For that purpose we denote by  $\BB_{R_0,\Omega} $  the set of all linear combinations of    $u_k(x) = \Phi_\Omega( x - k h_x )$ whose centers satisfy $\norm{k h_x} \leq R_0$. Moreover, we choose equidistant angular samples $\spos_m = R (\cos (m h_\theta ), \sin (m h_\theta ))$ for $m  =  0, \dots , M - 1$ for some angular step size $h_\theta  =   2 \pi / M $ .

\begin{theorem}[Angular sampling] \label{thm:angular}
Samples  $\Wo  f(\spos_m, \edot)$  stably determine all $f  \in \BB_{R_0,\Omega}$ if and only if $h_\theta  \leq \pi / (R_0\Omega)$.  
\end{theorem}

Theorems  \ref{thm:temporal}-\ref{thm:angular} complete the picture on  sampling in PAT and its  theoretical  analysis. Let an  initial pressure $f$ be of essential bandwidth $\Omega$ and located inside the disc of radius $R_0$.  Then the spatial step size  $h_x = \pi / \Omega$ allows to resolve $f$. Moreover,  temporal sampling rate $h_t = h_x$ and angular sampling rate $h_\theta = h_x / R_0$ are the minimal conditions that allow stable reconstruction of that function from sampled PAT data.  Taking angular samples is often costly or time consuming and therefore  undersampling is common, resulting in  sparse view PAT.

\begin{figure}[htb!]
\centering
\includegraphics[width=\textwidth]{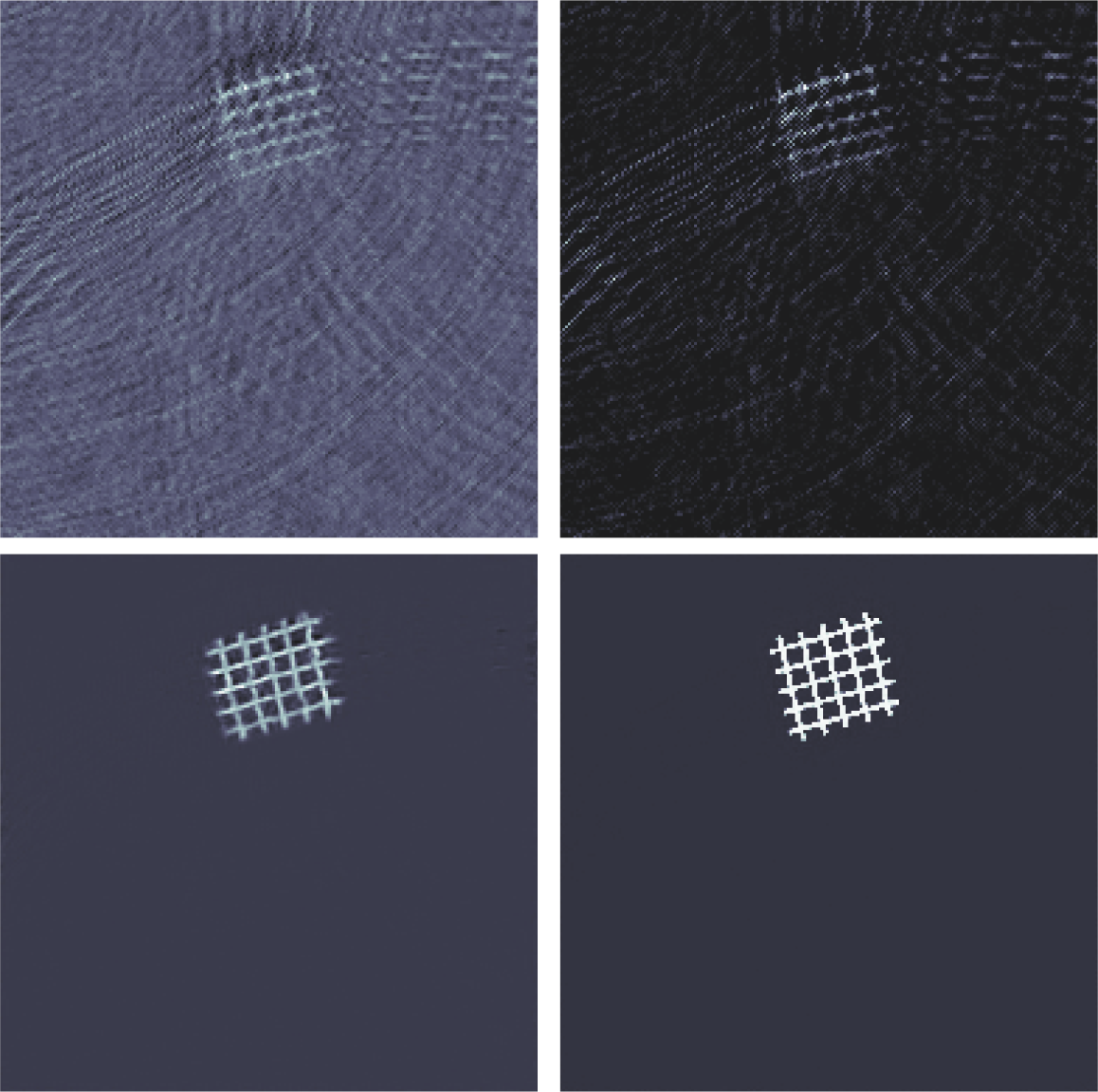}
\caption{Reconstructions from experimental data on a square of side length \SI{40}{\mm}. From left to right: Tikhonov regularization, joint $\ell^1$-minimization, learned primal dual iteration, and trained Unet.} \label{fig:results}
\end{figure}

\section{Experiment}

Experimental data have been acquired by an all-optical PAT projection imaging system 
described in \cite{bauer2017all}. The geometry of data acquisition and reconstruction is shown in Figure~\ref{fig:setup}.  The data is lowpass-filtered and downsampled. The system matrix is build  of the form $\Ao  =  \So \circ \Wo \circ \Uo^\ast$. Here  $\Uo^\ast$ is the synthesis operator for the translates $u_k = u(x - k h_x)$, where $h_x    \simeq   \SI{0.21}{\mm}  $ is the spatial step size and $k h_x$  is supported inside the square of side length \SI{40}{\mm}  centered at the origin.   This results in a total number of $192$ spatial samples in each coordinate direction. We use $ u(\rpos) =   ((\nu+1)/ \pi h_x^2) \,  (1- \norm{\rpos}^2/h_x^2)^\nu$ for $\norm{\rpos} \leq h_x$ and $u(\rpos) = 0$ otherwise with $\nu = 2$.  We use temporal sampling rate  (after temporal rescaling with the sound speed) $h_t = h_x$.  The  system uses 64 sensors arranged on the circle of radius  \SI{40}{\mm}  uniformly covering an angular range of   \SI{289}{\degree}. This results in an angular sampling step size  $h_\theta \simeq \SI{0.078}{\radian}$. The sampling condition is thus only satisfied inside the disc of radius   $h_x / h_\theta =  2.7$ which by far does not contain the initial pressure resulting in severe angular undersampling.   Resolving the  grid requires sampling rate  $h_{\rm grid} =   1/9 \si{\mm} $ that is satisfied for the angular sampling within a disc of radius  \SI{14}{\mm}.      Reconstruction results are shown in Figure \ref{fig:results}, where we compare standard quadratic Tikhonov regularization,  joint $\ell^1$-minimization   \cite{sandbichler2018sparsification}, learned primal dual  \cite{adler2018learned} and a trained Unet 
\cite{ronneberger2015u,antholzer2018deep}.

\section{Discussion}

All tested reconstruction methods are able to  recover the grid phantom  from angularly undersampled data. The quality of the machine learning methods is best. However no structures of the grid phantom seem to be lost even for standard quadratic  Tikhonov regularization. In the full proceedings we will make a detailed resolution study by further decreasing angular sampling and varying  locations of the grid phantom. In particular,  learned reconstruction methods will be critically analyzed  wether they provide improved resolution in a reliable and stable manner.

\section*{Acknowledgment}
This work has been supported by the Austrian Science Fund (FWF), project P 30747-N32.

\section*{Acknowledgements}
All authors acknowledge support of the Austrian Science Fund (FWF), project P 30747-N32.

\end{document}